\documentclass{amsart}

\usepackage{amsfonts,amssymb}

\newtheorem{corollary}{Corollary}
\newtheorem{theorem}{Theorem}
\newtheorem*{theorem*}{Theorem}
\newtheorem{lemma}{Lemma}

\newtheorem*{proposition*}{Proposition}

\newcommand{\rsl}{\mathfrak{sl}}

\newcommand{\sym}{\mathcal{H}}
\newcommand{\sks}{\mathcal{K}}
\newcommand{\wh}[1]{\widehat{#1}}
\newcommand{\wt}[1]{\widetilde{#1}}
\newcommand{\wb}[1]{\overline{#1}}
\newcommand{\vphi}{\varphi}
\newcommand{\iv}[1]{#1^{-1}\!}
\newcommand{\veps}{\varepsilon}
\newcommand{\chr}[1]{\mathrm{char}\,#1}
\newcommand{\ad}[1]{\mathrm{ad}\,#1}
\newcommand{\supp}[1]{\mathrm{Supp}\,#1}
\newcommand{\aut}[1]{\mathrm{Aut}\,#1}

\newcommand{\wg}{\widehat{G}}

\begin{document}

%
%Title and short title are preliminary
%
\title[Group gradings on simple Lie algebras]{Group gradings on simple Lie algebras of type A in positive characteristic}

\author[Bahturin]{Yuri Bahturin}
\address{Department of Mathematics and Statistics\\Memorial
University of Newfoundland\\ St. John's, NL, A1C5S7, Canada}
\email{yuri@math.mun.ca}

\author[Kochetov]{Mikhail Kochetov}
\address{Department of Mathematics and Statistics\\Memorial
University of Newfoundland\\ St. John's, NL, A1C5S7, Canada}
\email{mikhail@math.mun.ca}

\author[Montgomery]{Susan Montgomery}\address{Department of Mathematics
\\University of Southern California\\3620 South Vermont Ave., KAP 108\\
Los Angeles, CA 90089-2532, USA}\email{smontgom@math.usc.edu}

\thanks{The first author was partially supported by
NSERC grant \# 227060-04 and by an URP grant, Memorial University of
Newfoundland. The second author was supported by a Start-up grant, Memorial University of
Newfoundland. The third author was supported by NSF grant DMS 0401399.}

\begin{abstract}
In this paper we consider gradings by a finite abelian group $G$ on the Lie algebra $\rsl_n(F)$ over an algebraically
closed field $F$ of characteristic different from 2 and not dividing $n$. 
%The case of characteristic zero was dealt with in \cite{BZA}. 
\end{abstract}

%\begin{keyword}
%Graded algebra \sep matrix algebra\sep involution
%\end{keyword}

\maketitle

\section{Introduction}\label{si}

All gradings on the full matrix algebra $R=M_n(F)$ over an algebraically closed field $F$ by any finite group $G$ have been described in \cite{BSZ, BZnc, surgrad}. Namely, $R=A\otimes B$ where $A$ and $B$ are graded unital subalgebras such that $A\cong M_k(F)$ with a {\em fine} grading and $B\cong M_l(F)$ with an {\em elementary} grading defined by an $l$-tuple $(g_1,\ldots,g_l)$ so that if $a\in A_h$ then $a\otimes E_{ij}\in R_{\iv{g_i}hg_j}$. If $G$ is abelian, the case which is the most important for us, then $A$ decomposes as the tensor product of so-called $\veps$-graded matrix algebras, spanned by so-called {\em generalized Pauli matrices}. Note that the support $\supp{A}$ is always a subgroup $T$ of $G$. Also note that if $\chr{F}=p$ and $|G|$ or $n$ is a power of $p$, then all $G$-gradings of $M_n(F)$ are elementary.

We are interested in gradings on finite-dimensional simple Lie algebras. In the case $\chr{F}=0$, all gradings on the classical simple Lie algebras (except of type $D_4$) have been described in \cite{BShZ,BZA,antaut}. Here we will focus on the case $\chr{F}=p>0$. If a simple Lie algebra $L$ is graded by a group $G$, then $\supp{L}$ generates an abelian subgroup in $G$ (see \cite[Lemma 2.1]{BZA}). Thus it is sufficient to consider the case when $G$ is abelian. 

\section{Duality}\label{duality}

Let $G$ be a finite group, $F$ an algebraically closed  field. Let $H=FG$ be the group algebra of $G$ viewed as a Hopf algebra with comultiplication $\Delta(g)=g\otimes g$,  counit $\veps(g)=1$, and antipode $S(g)=g^{-1}$, for any $g\in G$. We will use Sweedler's notation: $\Delta(h)=\sum h_1\otimes h_2$, for any $h\in H$ (for the basic facts on Hopf algebras the reader is referred to \cite{Mont}).

Let $A$ be a nonassociative algebra over $F$. It is well-known that a $G$-grading on $A$ is equivalent to the structure of a right $H$-comodule algebra, i.e., a homomorphism of algebras $\rho: A\rightarrow A\otimes H$, written as $\rho(a)=\sum a_0\otimes a_1$ where $a_0\in A$ and $a_1\in H$, such that 
\[
(\rho\otimes id)\rho=(id\otimes \Delta)\rho.
\] 
Namely, if $A=\bigoplus_{g\in G}A_g$ is a $G$-graded algebra, then the mapping $\rho$ is defined on a homogeneous element $a$ of degree $g$ by $\rho(a)=a\otimes g$. 
Conversely, given $\rho: A\rightarrow A\otimes H$, one can define a $G$-grading on $A$ by setting $A_g=\{a\in A\,|\,\rho(a)=a\otimes g\}$, for any $g\in G$. 

Consider the dual space $K=H^*$. It has the natural structure of a Hopf algebra, with multiplication $(f' f'')(h)=\sum f'(h_1)f''(h_2)$, for all $h\in H$, and comultiplication 
$\Delta(f)=\sum f_1\otimes f_2$ if $\sum f_1(h')f_2(h'')=f(hk)$, for all $h',h''\in H$. 
In particular, if $\{e_g\,|\,g\in G\}$ is the basis of $K$ dual to $\{g\,|\,g\in G\}$, i.e., $e_g\in K$ are such that $e_g(h)=\delta_{g,h}$ for any $h\in G$ 
(Kronecker's delta), then $e_{g'} e_{g''}=\delta_{g',g''}e_{g'}$ and 
\[
\Delta(e_g)=\sum_{g',g''\in G:\, g'g''=g}e_{g'}\otimes e_{g''}.
\] 
It follows that $\sum_{g\in G}e_g$ is the unit element of $K$, $\veps(e_g)=\delta_{g,1}$, and $S(e_g)=e_{g^{-1}}$. Note that since $H$ is cocommutative, $K$ is commutative. In fact, $K$ is the direct product of fields $Fe_g$, $g\in G$. 

Now $K$ acts on $A$ by $f\cdot a=(id\otimes f)\rho(a)$. Using the definition of $\rho$, we obtain $f\cdot a=f(g)a$ for any $a\in A_g$, $g\in G$. With respect to this action $A$ becomes a $K$-module algebra, i.e.,  
\[
k\cdot(ab)
%=\mu(\Delta(k)(a\otimes b))
=\sum (k_1\cdot a)(k_2\cdot b)\mbox{ for all }k\in K,\; a,b\in A.
\]

Conversely, if $A$ is a $K$-module algebra, then there exists a homomorphism of algebras $\rho:A\rightarrow A\otimes H$ such that $K$ acts on $A$ by 
$f\cdot a=(id\otimes f)\rho(a)$. Also if $A$ is a unital (associative) algebra, then the requirement $1_A\in A_1$ is equivalent to $k\cdot 1_A=\veps(k) 1_A$ for all $k\in K$.

It is well-known that $B$ is a subcomodule of a right $H$-comodule $A$ if and only if $B$ is a submodule of the left $H^*$-module $A$. Here this fact means that $B$ is a graded subspace of $A$ if and only if $B$ is a $K$-submodule. 
%Indeed, if $B$ is a graded subspace then for any $b=\sum_{g\in G}b_g\in B$, where $b_g\in A_g$, we have $b_g\in B$ for all $g\in G$ and thus $\rho(b)=\sum_{g\in G}g\otimes b_g\in B\otimes H$, which means that $B$ is a subcomodule. Conversely, if $B$ is a subcomodule then using the same calculation and linear independence of $g\otimes b_g$ for different $g$ we conclude that each $b_g$ must be in $B$.

If $f\in K$ is a group-like element, i.e., $\Delta(f)=f\otimes f$ (hence $S(f)=f^{-1}$), then $f$ acts on $A$ as an automorphism: $f\cdot(ab)=(f\cdot a)(f\cdot b)$ for any $a,b\in A$.  The group-like elements of $K$ are the algebra homomorphisms $H\rightarrow F$, so their set can be identified with the group $\wg$ of multiplicative characters of $G$. If $\chi:G\rightarrow F^\times$ is a multiplicative character of $G$, then the element 
\[
\wt{\chi}=\sum_{g\in G} \chi(g)e_g
\]
is group-like in $K$. If $G$ has $|G|$ different characters, then the mapping $\chi\mapsto \wt{\chi}$ extends to an isomorphism of Hopf algebras $F\wh{G}\to K$. In this case $G$-gradings on an algebra $A$ are equivalent to $\wh{G}$-actions on $A$ by automorphisms. We have this situation if and only if $G$ is abelian and $\chr{F}$ does not divide $|G|$.

Now if $f\in K$ is primitive, i.e., $\Delta(f)=f\otimes 1+1\otimes f$, then $f$ acts on $A$ as a derivation: $f\cdot(ab)=(f\cdot a)b+a(f\cdot b)$ for any $a,b\in A$.
Let $\alpha:G\rightarrow F$ be a map and set 
\[
\wt{\alpha}=\sum_{g\in G}\alpha(g)e_g.
\]
It is easy to check that $\wt{\alpha}$ is primitive if and only if $\alpha:G\rightarrow F$ is an additive character of $G$. Of course, nonzero additive characters can 
exist only if $\chr{F}=p>0$ and $p$ divides $|G|$.

For example, let $G=\langle a_1\rangle_p\times\cdots\times\langle a_k\rangle_p$, an elementary abelian $p$-group. Then there exist $k$ additive characters  $\alpha_1,\ldots,\alpha_k$ such that $\alpha_i(a_j)=\delta_{i,j}$. Then the elements $\wt{\alpha}_i$ are primitive and also satisfy $(\wt{\alpha}_i)^p=\wt{\alpha}_i$. 
The span of the elements $\wt{\alpha}_i$ in $K$ is an abelian $p$-Lie algebra $\mathfrak{g}$, and $K$ is isomorphic to the restricted enveloping algebra $u(\mathfrak{g})$. %Indeed, let us consider the standard basis of $u(\mathfrak{g})$ which consists of all monomials $\wt{\alpha}_1^{m_1}\cdots\wt{\alpha}_k^{m_k}$, $0\leq m_i < p$. The result of the action of this monomial on $g=a_1^{l_1}\cdots a_k^{l_k}$ can be easily computed to be $l_1^{m_1}\cdots l_k^{m_k}$. Thus a linear combination of such monomials is zero if and only if a polynomial in $k$ variables of degree less than $p$ in each variable in the polynomial ring $F[x_1,\ld,x_k]$ vanishes in all integral points $(l_1,\ldots,l_k)\in \mathbb{Z}_p^k$. Using Vandermonde's determinants one easily sees that all coefficient of such polynomial are zero. This in its turn proves that all basis monomials of $u(\mathfrak{g})$ are linearly independent in $K$ and so the elements $\wt{\alpha}_i$, $i=1,\ld,k$ generate in $K$ the subalgebra isomorphic to $u(\mathfrak{g})$.
In this case $G$-gradings on an algebra $A$ are equivalent to $\mathfrak{g}$-actions on $A$ by derivations.

Now let $G$ be any finite abelian group and $F$ an algebraically closed field of characteristic $p>0$. Then $H=FG$ and $K=H^*$ are finite-dimensional commutative and cocommutative Hopf algebras. We can write $G=G_0\times G_1$ where $G_0$ is of order not divisible by $p$ and $G_1$ is a $p$-group. This induces the following decompositions of $H$ and $K$: 
$H=H_0\otimes H_1$ where $H_0=FG_0$ and $H_1=FG_1$, and $K=K_0\otimes K_1$ where $K_0=(H_0)^*$ and $K_1=(H_1)^*$. By duality, a $G$-grading on $L$ is equivalent to the structure of a $K$-module algebra. Since $K=K_0\otimes K_1$, the latter is equivalent to a pair of mutually commuting actions on $L$ by $K_0$ and by $K_1$ 
that make $L$ a $K_0$-module algebra, resp., $K_1$-module algebra. The $K_0$-module structure on $L$ is equivalent to a $\wh{G}_0$-action on $L$ by automorphisms. If $G_1$ is an elementary abelian $p$-group, then the $K_1$-module structure on $L$ is equivalent to a $\mathfrak{g}_1$-action by derivations, where $\mathfrak{g}_1$ is the abelian $p$-Lie algebra associated to $G_1$. If $G_1$ is not elementary, the situation is more complicated and involves the so-called divided power algebras, which will be discussed in the next section.

\section{Hopf actions on matrix algebras}\label{Hopf}

Let $R=M_n(F)$ where $F$ is a field of characteristic $p>0$. Let $G$ be a finite abelian $p$-group. We want to describe all $G$-gradings on the Lie algebra $R^{(-)}$. 
Set $H=FG$ and $K=H^*$. As discussed in the previous section, a $G$-grading on $R$, resp. $R^{(-)}$, is equivalent to a $K$-comodule algebra structure on $R$, resp. $R^{(-)}$.
If $G$ is an elementary abelian $p$-group of rank $k$, then $K=u(\mathfrak{g})$ where $\mathfrak{g}$ is an abelian $p$-Lie algebra of dimension $k$. Any element 
$\delta\in\mathfrak{g}$ then acts as a derivation of $R$, resp. $R^{(-)}$. Thus we can apply a result of Martindale \cite{Mart} on Lie derivations of a primitive ring with a 
nontrivial idempotent. Here we need the result only in the case of a simple ring:

\begin{theorem*}[Martindale]
Let $R$ be a simple associative unital ring. Assume that the characteristic of $R$ is not $2$ and $R$ contains a nontrivial idempotent. Let $\delta:R\rightarrow R$ be a derivation of $R^{(-)}$. Then $\delta=\tau+\zeta$ where $\tau:R\rightarrow R$ is a derivation of $R$ and $\zeta:R\rightarrow Z(R)$ is an additive map that vanishes on $[R,R]$.
\end{theorem*}

Applying the above theorem to $R=M_n(F)$ ($n\geq 2$), we obtain the following:

\begin{corollary}\label{c1}
Let $R=M_n(F)$, $\chr{F}=p>0$, $p\neq 2$ and $p\nmid n$. Let $G$ be an elementary abelian $p$-group. Suppose $R=\bigoplus_{g\in G}R_g$ is a grading on $R^{(-)}$. 
Then $R=\bigoplus_{g\in G}R_g$ is a grading on $R$ if and only if $1\in R_1$.
\end{corollary}

\begin{proof} By the above discussion, we have that the abelian $p$-Lie algebra $\mathfrak{g}$, corresponding to the group $G$, acts on $R$ by Lie derivations. Moreover, each of these Lie derivations maps $1$ to $0$. It follows by the Martindale's theorem that these operators are in fact associative derivations of $R$. Obviously, they continue to satisfy the same relations when regarded as associative derivations. Therefore, the associative algebra $R$ is a $u(\mathfrak{g})$-module algebra and thus a $G$-graded algebra (with the same subspaces $R_g$ as for $R^{(-)}$). \end{proof}

We want to extend this result to an arbitrary finite abelian $p$-group $G$. First consider the case $G=\langle a \rangle_{p^N}$. Then $H=F[t]/(t^{p^N}-1)=F[\xi]/(\xi^{p^N})$ where $\xi=t-1$. Let $\{\delta^{(m)}\,|\,m=0,\ldots,p^N-1\}$ be the basis of $K$ dual to $\{\xi^m\,|\,m=0,\ldots,p^N-1\}$. Then the coproduct of $K$ is given by
\begin{equation}\label{div_powers}
\Delta\delta^{(m)}=\sum_{i=0}^{m}\delta^{(i)}\otimes\delta^{(m-i)}.
\end{equation}
Elements $\delta^{(m)}$ with coproduct of this form are sometimes called ``divided powers''. In particular, $\delta^{(0)}=1$ and $\delta^{(1)}$ spans the space of primitive elements of $K$. One can also write an explicit formula for the product $\delta^{(i)}\delta^{(j)}$, but we will only need that
\[
\delta^{(i)}\delta^{(j)}=\binom{i+j}{i}\delta^{(i+j)}\pmod{{\rm span}\{\delta^{(m)}\,|\,m<i+j\}}.
\]
It follows that, for any $1\leq l\leq N$, the subspace spanned by $\delta^{(m)}$ with $m<p^l$ is a subalgebra of $K$, which is generated by the elements $\delta^{(p^k)}$, $k=0,\ldots,l-1$ (see e.g. \cite{Dieu}). In particular, for $N>1$ the algebra $K$ is not generated by primitive elements and, consequently, we will have to consider operators 
with more complicated ``product expansion laws'' than the ordinary Leibniz rule (see e.g. (\ref{gen_Lie_der}) below).

\begin{theorem}\label{Hopf_action}
Let $F$ be a field of characteristic $p\neq 2$ and $R=M_n(F)$ with $p\nmid n$. Let $G$ be a finite abelian $p$-group and $K=(FG)^*$. 
Then any $K$-module algebra structure on $R^{(-)}$ with the property that $1\in R$ is $K$-invariant, is in fact a $K$-module algebra structure on $R$.
\end{theorem}

Using duality, we can immediately reformulate the above theorem as follows (with the addition of an obvious ``only if'' part): 

\begin{corollary}\label{p_grading}
Let $R=M_n(F)$, $\chr{F}=p>0$, $p\neq 2$ and $p\nmid n$. Let $G$ be a finite abelian $p$-group. Suppose $R=\bigoplus_{g\in G}R_g$ is a grading on $R^{(-)}$. 
Then $R=\bigoplus_{g\in G}R_g$ is a grading on $R$ if and only if $1\in R_1$. $\hfill{\square}$
\end{corollary}

\begin{proof}
We will proceed by induction on $|G|$. We start by separating one cyclic factor: $G=\langle a \rangle_{p^N}\times\wt{G}$, hence $H=F\langle a \rangle\otimes\wt{H}$ and 
$K=(F\langle a \rangle)^*\otimes\wt{K}$. We introduce $\delta^{(m)}$ in the first factor as discussed above for the case of a cyclic group. Let $\wb{K}$ be generated by $\wt{K}$ and $\delta^{(p^k)}$, $k=0,\ldots,N-2$. Then $\wb{K}=(F\wb{G})^*$ where $\wb{G}=G/\langle a^p\rangle$ is a group of smaller order. 

By inductive hypothesis, $R$ is a $\wb{K}$-module algebra. This means that the factor-grading by $\wb{G}$: 
$R=\bigoplus_{\bar{g}\in\wb{G}}R_{\bar{g}}$ where $R_{\bar{g}}=\bigoplus_{g\in\bar{g}}R_g$, is a grading of $R$ as an associative algebra and hence an elementary grading. 
It follows that $R_{\bar{1}}$ contains a nontrivial idempotent $e$ (we assume $n\geq 2$, the case of $n=1$ being trivial). In the dual language, $e$ is $\wb{K}$-invariant, i.e., 
\begin{equation}\label{invar_idemp}
\delta^{(m)}\cdot e=0\quad\mbox{for}\quad m=1,\ldots,p^{N-1}-1.
\end{equation} 

Let $q=p^{N-1}$ and consider the operator $\sigma:R\rightarrow R$ defined by $\sigma(x)=\delta^{(q)}\cdot x$. Since $R^{(-)}$ is a $K$-module algebra, 
(\ref{div_powers}) implies 
\begin{equation}\label{gen_Lie_der}
\sigma([x,y])=[\sigma(x),y]+[x,\sigma(y)]+\sum_{k=1}^{q-1} [\delta^{(k)}\cdot x,\delta^{(q-k)}\cdot y]\qquad\forall x,y\in R
\end{equation}
and similarly for three or more factors. The goal is to show that 
\begin{equation}\label{gen_der}
\sigma(xy)=\sigma(x)y+x\sigma(y)+\sum_{k=1}^{q-1} (\delta^{(k)}\cdot x)(\delta^{(q-k)}\cdot y)\qquad\forall x,y\in R,
\end{equation}
which will mean that $R$ is a $K$-module algebra. Note that we already know that the analogue of (\ref{gen_der}) holds for $\delta^{(m)}$, $m<q$:
\begin{equation}\label{ind_hyp}
\delta^{(m)}\cdot(xy)=(\delta^{(m)}\cdot x)y+x(\delta^{(m)}\cdot y)+\sum_{k=1}^{m-1} (\delta^{(k)}\cdot x)(\delta^{(m-k)}\cdot y)\qquad\forall x,y\in R,
\end{equation}
since $R$ is a $\wb{K}$-module algebra.

The proof consists of a sequence of lemmas that are adaptations of those found in \cite{Mart} to our situation.
We will use the following notation. Set $e_1=e$, $e_2=1-e$, so $R=\bigoplus_{i,j=1}^2 R_{ij}$ where $R_{ij}=e_i R e_j$.

\begin{lemma}\label{lemma1}
$\sigma(e)=[s,e]+z$ for some $s\in R$ and $z\in Z(R)$.
\end{lemma}

\begin{proof}
Since $e^2=e$, we have $[e,[e,[e,x]]]=[e,x]$ for all $x\in R$. If we apply $\sigma$ to both sides and expand using (\ref{gen_Lie_der}) and taking into account 
(\ref{invar_idemp}) the we obtain the following.
\[
[\sigma(e),[e,[e,x]]]+[e,[\sigma(e),[e,x]]]+[e,[e,[\sigma(e),x]]]+[e,[e,[e,\sigma(x)]]]=[\sigma(e),x]+[e,\sigma(x)].
\]
Cancelling $[e,[e,[e,\sigma(x)]]]$ with $[e,\sigma(x)]$ and expanding the commutators, we obtain:
\begin{eqnarray}\label{sigma_idemp}
&&(e\sigma(e)+\sigma(e)e+e\sigma(e)e-\sigma(e))x-x(e\sigma(e)+\sigma(e)e+e\sigma(e)e-\sigma(e))\\
&&=3(e\sigma(e)+\sigma(e)e-\sigma(e))xe-3ex(e\sigma(e)+\sigma(e)e-\sigma(e))\qquad\forall x\in R.\nonumber
\end{eqnarray}
Write $\sigma(e)=\sum f_{ij}$ where $f_{ij}\in R_{ij}$.
Substituting into (\ref{sigma_idemp}) and simplifying, we obtain:
\begin{equation}\label{f_ij}
(2f_{11}-f_{22})x-x(2f_{11}-f_{22})=3(f_{11}-f_{22})xe-3ex(f_{11}-f_{22}).
\end{equation}

For $x\in R_{12}$, (\ref{f_ij}) gives $2f_{11}x+xf_{22}=3xf_{22}$, hence $f_{11}x=xf_{22}$ ($\chr{F}\neq 2$). Similarly, for $x\in R_{21}$, we get $f_{22}x=xf_{11}$. 
Set $z=f_{11}+f_{22}$. Then $z$ commutes with $R_{12}$ and $R_{21}$.

Now fix $x\in R_{11}$. For any $y\in R_{12}$, we have $(zx-xz)y=z(xy)-x(zy)=(xy)z-x(yz)=0$. Since the left annihilator of $R_{12}$ in $R$ is $Re_2$ and $zx-xz\in Re_1$, 
we conclude that $zx-xz=0$. Similarly, $zx-xz=0$ for $x\in R_{22}$.

We have proved that $z\in Z(R)$. Set $s=f_{21}-f_{12}$. Then $\sigma(e)=[s,e]+z$, as desired.
\end{proof}

Set $\wt{\sigma}=\sigma-\ad s$. Then $\wt{\sigma}(e)=z\in Z(R)$. Since $\ad s$ is a derivation of $R^{(-)}$, equation (\ref{gen_Lie_der}) holds with $\sigma$ replaced by $\wt{\sigma}$.

\begin{lemma}\label{lemma2}
$\wt{\sigma}(R_{ij})\subset R_{ij}$ for $i\neq j$.
\end{lemma}

\begin{proof}
We will show that $\wt{\sigma}(R_{12})\subset R_{12}$, the case of $R_{21}$ being similar. For any $x\in R_{12}$ we have $x=[e,x]$. Write $\wt{\sigma}(x)=\sum y_{ij}$ where $y_{ij}\in R_{ij}$. Then $\sum y_{ij}=\wt{\sigma}(x)=\wt{\sigma}([e,x])
=[\wt{\sigma}(e),x]+[e,\wt{\sigma}(x)]=[e,\wt{\sigma}(x)]$, where we used (\ref{gen_Lie_der}), (\ref{invar_idemp}), and 
$\wt{\sigma}(e)\in Z(R)$. Now $\sum y_{ij}=[e,\wt{\sigma}(x)]=[e,\sum y_{ij}]=y_{12}+y_{11}-(y_{11}+y_{21})=y_{12}-y_{21}$. 
It follows that $y_{11}=y_{22}=0$ and also $y_{21}=0$ ($\chr{F}\neq 2$).
\end{proof}

\begin{lemma}\label{lemma3}
$\wt{\sigma}(R_{ii})\subset R_{ii}\oplus Z(R)$.
\end{lemma}

\begin{proof}
Let $x\in R_{11}$. Write $\wt{\sigma}(x)=\sum y_{ij}$ where $y_{ij}\in R_{ij}$. Since $[e,x]=0$, we obtain
$0=\wt{\sigma}([e,x])=[e,\wt{\sigma}(x)]=y_{12}-y_{21}$ and hence $y_{12}=y_{21}=0$. 
Thus $\wt{\sigma}(R_{11})\subset R_{11}\oplus R_{22}$. Similarly, $\wt{\sigma}(R_{22})\subset R_{11}\oplus R_{22}$.

Fix $x\in R_{11}$ and $y\in R_{22}$. Since $[x,y]=0$, we have 
\[
[\wt{\sigma}(x),y]+[x,\wt{\sigma}(y)]+\sum_{k=1}^{q-1}[\delta^{(k)}\cdot x,\delta^{(q-k)}\cdot y]=0.
\]
From (\ref{invar_idemp}) and the fact that $R$ is a $\wb{K}$-module algebra it follows that each $R_{ij}$ is 
$\delta^{(k)}$-invariant for $i,j=1,2$ and $k=1,\ldots,q-1$. Hence $[\delta^{(k)}\cdot x,\delta^{(q-k)}\cdot y]=0$ 
for all $k$ and we obtain
\[
[\wt{\sigma}(x),y]+[x,\wt{\sigma}(y)]=0.
\]

Now write $\wt{\sigma}(x)=a_{11}+a_{22}$ and $\wt{\sigma}(y)=b_{11}+b_{22}$ where $a_{11},b_{11}\in R_{11}$ and $a_{22},b_{22}\in R_{22}$. Then the above equation gives
\[
[a_{22},y]+[x,b_{11}]=0
\]
where the first term is in $R_{22}$ and the second is in $R_{11}$, so $[a_{22},y]=0$ and $[x,b_{11}]=0$.

We have proved that $[a_{22},y]=0$ for all $y\in R_{22}$. Therefore, $a_{22}=\lambda(1-e)$ for some $\lambda\in F$. 
Hence $\wt{\sigma}(x)=a_{11}+\lambda(1-e)=a'_{11}+\lambda 1$ where $a'_{11}=a_{11}-\lambda e\in R_{11}$. 
Similarly, $\wt{\sigma}(y)=b'_{22}+\mu 1$ where $b'_{22}\in R_{22}$ and $\mu\in F$.
\end{proof}

From Lemmas \ref{lemma2} and \ref{lemma3} (and the fact that $R_{ii}\cap Z(R)=0$) it follows that we can uniquely define 
linear maps $\tau$ and $\zeta$ on $R$ such that $\wt{\sigma}=\tau+\zeta$, $\tau(R_{ij})\subset R_{ij}$ for $i,j=1,2$, and $\zeta(R)\subset Z(R)$. From Lemma \ref{lemma2} we also have $\zeta(R_{ij})=0$ for $i\neq j$.

\begin{lemma}\label{lemma4}
For all $x\in R_{ij}$ with $i\neq j$ and all $y\in R$,
\[
\tau(xyx)=\tau(x)yx+x\tau(y)x+xy\tau(x)+\sum_{\stackrel{0\leq k,l,m < q}{k+l+m=q}}
(\delta^{(k)}\cdot x)(\delta^{(l)}\cdot y)(\delta^{(m)}\cdot x).
\]
\end{lemma}

\begin{proof}
Observe that $R^2_{ij}=0$ implies that $x_1 y x_2 + x_2 y x_1=[x_1,[y,x_2]]$ for all $x_1,x_2\in R_{ij}$ and $y\in R$.
Using this and the analogue of (\ref{gen_Lie_der}) for three factors, we compute:
\begin{eqnarray*}
&&2\tau(xyx)=\wt{\sigma}(2xyx)=\wt{\sigma}([x,[y,x]])\\
&=&[\wt{\sigma}(x),[y,x]]+[x,[\wt{\sigma}(y),x]]+[x,[y,\wt{\sigma}(x)]]
+\sum_{\stackrel{0\leq k,l,m < q}{k+l+m=q}}[\delta^{(k)}\cdot x,[\delta^{(l)}\cdot y,\delta^{(m)}\cdot x]]\\
&=&[\tau(x),[y,x]]+[x,[\tau(y),x]]+[x,[y,\tau(x)]]
+\sum_{\stackrel{0\leq k,l,m < q}{k+l+m=q}}[\delta^{(k)}\cdot x,[\delta^{(l)}\cdot y,\delta^{(m)}\cdot x]]\\
&=&\tau(x)yx+xy\tau(x)+2x\tau(y)x+xy\tau(x)+\tau(x)yx\\
&&+\sum_{\stackrel{0\leq k,l,m < q}{k+l+m=q}}(\delta^{(k)}\cdot x)(\delta^{(l)}\cdot y)(\delta^{(m)}\cdot x)
+\sum_{\stackrel{0\leq k,l,m < q}{k+l+m=q}}(\delta^{(m)}\cdot x)(\delta^{(l)}\cdot y)(\delta^{(k)}\cdot x)\\
&=&2\tau(x)yx+2x\tau(y)x+2xy\tau(x)
+2\sum_{\stackrel{0\leq k,l,m < q}{k+l+m=q}}(\delta^{(k)}\cdot x)(\delta^{(l)}\cdot y)(\delta^{(m)}\cdot x).
\end{eqnarray*}
It remains to cancel 2.
\end{proof}

\begin{lemma}\label{lemma5}
If $x\in R_{ii}$, $y\in R_{jk}$, and $j\neq k$, then 
\begin{equation}\label{gen_der_}
\tau(xy)=\tau(x)y+x\tau(y)+\sum_{k=1}^{q-1} (\delta^{(k)}\cdot x)(\delta^{(q-k)}\cdot y).
\end{equation}
\end{lemma}

\begin{proof}
If $i\neq j$, then both sides are zero, so without loss of generality we assume $x\in R_{11}$ and $y\in R_{12}$.
Then $xy=[x,y]$ and hence
\begin{eqnarray*}
&&\tau(xy)=\wt{\sigma}(xy)=\wt{\sigma}([x,y])\\
&&=[\wt{\sigma}(x),y]+[x,\wt{\sigma}(y)]+\sum_{k=1}^{q-1} [\delta^{(k)}\cdot x,\delta^{(q-k)}\cdot y]\\
&&=[\tau(x),y]+[x,\tau(y)]+\sum_{k=1}^{q-1} [\delta^{(k)}\cdot x,\delta^{(q-k)}\cdot y]\\
&&=\tau(x)y+x\tau(y)+\sum_{k=1}^{q-1} (\delta^{(k)}\cdot x)(\delta^{(q-k)}\cdot y).
\end{eqnarray*}
\end{proof}

\begin{lemma}\label{lemma6}
If $x\in R_{ii}$ and $y\in R_{jj}$, then equation (\ref{gen_der_}) holds.
\end{lemma}

\begin{proof}
If $i\neq j$, then both sides are zero, so we assume without loss of generality $x,y\in R_{11}$. Fix any $r\in R_{12}$.
Then from Lemma \ref{lemma5} we obtain
\begin{eqnarray*}
\tau(xy)r&=&\tau(xyr)-xy\tau(r)-\sum_{l=1}^{q-1} (\delta^{(l)}\cdot (xy))(\delta^{(q-l)}\cdot r)\\
&=&\tau(x)yr+x\tau(yr)+\sum_{k=1}^{q-1} (\delta^{(k)}\cdot x)(\delta^{(q-k)}\cdot (yr))\\
&&-xy\tau(r)-\sum_{l=1}^{q-1} (\delta^{(l)}\cdot (xy))(\delta^{(q-l)}\cdot r)\\
&=&\tau(x)yr+x\tau(y)r+xy\tau(r)+x\sum_{m=1}^{q-1} (\delta^{(m)}\cdot y)(\delta^{(q-m)}\cdot r)\\
&&+\sum_{k=1}^{q-1} (\delta^{(k)}\cdot x)(\delta^{(q-k)}\cdot (yr))
-xy\tau(r)-\sum_{l=1}^{q-1} (\delta^{(l)}\cdot (xy))(\delta^{(q-l)}\cdot r).
\end{eqnarray*}
Now we simplify our expression by expanding the products in the last two summations using (\ref{ind_hyp}) 
and then cancelling the common terms. This yields
\begin{eqnarray*}
\tau(xy)r&=&(\tau(x)y+x\tau(y))r+x\sum_{m=1}^{q-1} (\delta^{(m)}\cdot y)(\delta^{(q-m)}\cdot r)\\
&&+\sum_{k=1}^{q-1} (\delta^{(k)}\cdot x)(\delta^{(q-k)}\cdot y)r
-\sum_{l=1}^{q-1} x(\delta^{(l)}\cdot y)(\delta^{(q-l)}\cdot r)\\
&=&\left(\tau(x)y+x\tau(y)+\sum_{k=1}^{q-1} (\delta^{(k)}\cdot x)(\delta^{(q-k)}\cdot y)\right)r.
\end{eqnarray*}
Since $r\in R_{12}$ was arbitrary, we obtain the desired equation.
\end{proof}

Now we will show that equation (\ref{gen_der_}) holds for arbitrary $x,y\in R$. Lemmas \ref{lemma5} and \ref{lemma6}
(and also an analogue of Lemma \ref{lemma5} where $x$ and $y$ are interchanged) cover all cases except the following: 
$x\in R_{ij}$ and $y\in R_{st}$ where $i\neq j$ and $s\neq t$. If $j\neq s$, then both sides of (\ref{gen_der_}) are zero, 
so it remains to consider only two possibilities: $x\in R_{12}$, $y\in R_{21}$ or $x\in R_{21}$, $y\in R_{12}$. 
Interchanging $x$ and $y$, it suffices to consider only the first possibility: $x\in R_{12}$, $y\in R_{21}$. 
Set $z=\zeta([x,y])$. Then 
\begin{eqnarray*}
&&z=\wt{\sigma}([x,y])-\tau([x,y])\\
&&=[\wt{\sigma}(x),y]+[x,\wt{\sigma}(y)]+\sum_{l=1}^{q-1} [\delta^{(l)}\cdot x,\delta^{(q-l)}\cdot y]-\tau([x,y])\\
&&=[\tau(x),y]+[x,\tau(y)]-\tau([x,y])+\sum_{l=1}^{q-1} [\delta^{(l)}\cdot x,\delta^{(q-l)}\cdot y]
\end{eqnarray*}
Expanding the commutators and grouping the terms, we obtain:
\begin{eqnarray}\label{two_components}
z&=&\left(\tau(x)y+x\tau(y)-\tau(xy)+\sum_{l=1}^{q-1} (\delta^{(l)}\cdot x)(\delta^{(q-l)}\cdot y)\right)\\
&&-\left(\tau(y)x+y\tau(x)-\tau(yx)+\sum_{l=1}^{q-1} (\delta^{(l)}\cdot y)(\delta^{(q-l)}\cdot x)\right)\nonumber
\end{eqnarray}
The first group of terms is in $R_{11}$ and the second is in $R_{22}$. So if we prove that $z=0$, then equation (\ref{gen_der_}) will follow.

Multiplying both sides of (\ref{two_components}) by $x$ on the left, we obtain:
\[
xz=x\tau(yx)-x\tau(y)x-xy\tau(x)-x\sum_{l=1}^{q-1} (\delta^{(l)}\cdot y)(\delta^{(q-l)}\cdot x). 
\]
By Lemma \ref{lemma5}, 
\[
x\tau(yx)=\tau(xyx)-\tau(x)yx-\sum_{k=1}^{q-1} (\delta^{(k)}\cdot x)(\delta^{(q-k)}\cdot (yx)).
\]
Therefore,
\begin{eqnarray*}
xz&=&\tau(xyx)-\tau(x)yx-x\tau(y)x-xy\tau(x)\\
&&-\sum_{l=1}^{q-1} x(\delta^{(l)}\cdot y)(\delta^{(q-l)}\cdot x)
-\sum_{k=1}^{q-1} (\delta^{(k)}\cdot x)(\delta^{(q-k)}\cdot (yx))\\
&=&\tau(xyx)-\tau(x)yx-x\tau(y)x-xy\tau(x)
-\sum_{\stackrel{0\leq k,l,m < q}{k+l+m=q}}(\delta^{(k)}\cdot x)(\delta^{(l)}\cdot y)(\delta^{(m)}\cdot x).
\end{eqnarray*}
By Lemma \ref{lemma4}, we conclude that $xz=0$. But $z=\lambda 1$ for some $\lambda\in F$, so $z\neq 0$ would imply 
$x=0$ and hence $z=\zeta([x,y])=0$ --- a contradiction. Therefore, $z=0$ and we have proved equation (\ref{gen_der_}) for all $x,y\in R$.

It follows that $\zeta([x,y])=0$ for all $x,y\in R$. Since also $\zeta(1)=0$ and $p\nmid n$, we conclude that $\zeta=0$. 
Thus $\wt{\sigma}=\tau$ and equation (\ref{gen_der_}) reads
\[
\wt{\sigma}(xy)=\wt{\sigma}(x)y+x\wt{\sigma}(y)+\sum_{k=1}^{q-1} (\delta^{(k)}\cdot x)(\delta^{(q-k)}\cdot y).
\]
Recalling that $\sigma=\wt{\sigma}+\ad s$, we obtain (\ref{gen_der}), as desired. 

\end{proof}

\section{Gradings on $\rsl_n(F)$}

The gradings on the Lie algebra $L=\rsl_n(F)$ over an algebraically closed field $F$ of characteristic zero have been completely described in \cite{BZA}. 
Namely, the gradings $L=\bigoplus_{g\in G}L_g$ by a finite (abelian) group $G$ are of the following two types: 
\begin{enumerate}
\item[I\;:] $L_g=R_g$ for $g\neq 1$ and $L_1=R_1\cap L$ where $M_n(F)=\bigoplus_{g\in G}R_g$ is a $G$-grading on $M_n(F)$;
\item[II\;:] $L_g=\sks(R_g,*)\oplus \sym(R_{gh},*)$ if $g\neq h$ and $L_h=\sks(R_h,*)\oplus(\sym(R_1,*)\cap L)$ where $M_n(F)=\bigoplus_{g\in G}R_g$ is a $G$-grading on $M_n(F)$, $*$ is an involution that preserves the grading, and $h\in G$ is an element of order 2.
\end{enumerate}

The proof is based on the following key ideas. First, in this case the gradings by a finite abelian group $G$ are equivalent to the actions of $\wh{G}$ by automorphisms. Second, any inner automorphism of $\rsl_n(F)$ uniquely extends to an automorphism of $M_n(F)$ and any outer automorphism to the negative of an antiautomorphism of $M_n(F)$. Third, the antiautomorphisms of $M_n(F)$ that may arise here can be ``corrected'' by slightly changing the $\wh{G}$-action so they become automorphisms (see the proposition below). Finally, the original grading on $\rsl_n(F)$ can be recovered from the grading associated to the modified action on $M_n(F)$ by using an ``exchange formula'' 
(see Lemma \ref{exchange} below).

The goal of this section is to extend the above approach to describe the gradings on $\rsl_n(F)$ where $F$ is of positive characteristic $p\neq 2$ not dividing $n$. It turns out that in this case the answer is the same as in characteristic zero:

\begin{theorem}\label{main_sl} 
Let $L=\rsl_n(F)$ where $F$ is an algebraically closed field, $\chr{F}\neq 2$ and $\chr{F}\nmid n$. Let $G$ be a finite abelian group. Then any $G$-grading on $L$ is either of type I or of type II above. Moreover, if $G$ is a $p$-group then any $G$-grading on $L$ is of type I, i.e., the restriction of an \emph{elementary} $G$-grading of $M_n(F)$. 
\end{theorem}

\begin{proof}

As discussed in Section \ref{duality}, the gradings by $G$ are equivalent to the actions of the Hopf algebra $K=(FG)^*$. We write $G=G_0\times G_1$ where $G_0$ is of order not divisible by $p$ and $G_1$ is a $p$-group. Then $K=K_0\otimes K_1$ where $K_0=(FG_0)^*=F\wh{G_0}$ and $K_1=(FG_1)^*$. As in the case of characteristic zero, the action of $\wh{G_0}$ on $\rsl_n(F)$ can be extended to $M_n(F)$ thanks to the results of Blau and Martindale, summarized in \cite[Theorem 6.1]{BlMa} as follows:

\begin{theorem*}[Blau--Martindale] 
Let $S=M_m(E)$, $R=M_n(F)$, $n>1$, $E$ and $F$ fields with isomorphism $\gamma: F\rightarrow E$. Assume that $\chr{E}\neq 2$, and $m\neq 3$ if $\chr{E}=3$. Suppose there is a $\gamma$-semilinear Lie isomorphism $\alpha:\overline{[R,R]}\rightarrow \overline{[S,S]}$ where $\overline{[R,R]}=[R,R]/[R,R]\cap F$ and $\overline{[S,S]}=[S,S]/[S,S]\cap E$. Then $n=m$ and there exists a $\gamma$-semilinear map $\sigma:R\rightarrow S$ such that $\sigma$ is either an isomorphism or the negative of an antiisomorphism and such that $\overline{x^\alpha}=\overline{x}^\alpha$ for all $x\in [R,R]$.
\end{theorem*}

In our case, $E=F$, $\gamma=\mathrm{id}$, $R=S$, and $\chr{F}\nmid n$, so $\overline{[R,R]}=\rsl_n(F)$. Thus we can extend (uniquely) the action of $\wh{G_0}$ on $\rsl_n(F)$ to $R=M_n(F)$ and obtain a homomorphism $f:\wh{G_0}\to GL(R)$ whose image consists of automorphisms and, possibly, the negatives of antiautomorphisms of $R$, which are all 
automorphisms of the Lie algebra $R^{(-)}$. 

We also extend the $K_1$-action on $\rsl_n(F)$ to an action on $R^{(-)}$ by declaring that the identity matrix is $K_1$-invariant. Then by Theorem \ref{Hopf_action}, this action turns the associative algebra $R$ into a $K_1$-module algebra.

The extended action of $K_0\otimes K_1$ on $R^{(-)}$ corresponds to a Lie grading on $R$, $R=\bigoplus_{g\in G}R_g$, which restricts to the original $G$-grading on $\rsl_n(F)$. 

Now set $\Lambda=f^{-1}(\aut{R})$. This is a subgroup in $\wh{G_0}$ of index at most 2 that acts by automorphisms on $R$. Set $H=\Lambda^\perp$ in $G_0$. 
Then $H=\langle h \rangle$ where $h\in G_0$ is of order at most 2. Let $\wb{K}=F\Lambda\otimes K_1$. By construction, $R$ is a $\wb{K}$-module algebra, so the corresponding 
factor-grading by $\wb{G}=G/H$ on $R$ is a grading of $R$ as an associative algebra.

If $\Lambda=\wh{G_0}$, then we are done: we have a type I grading on $\rsl_n(F)$. Otherwise $\wh{G_0}$ is generated over $\Lambda$ by an element $\chi$ such that $f(\chi)=-\vphi$ where $\vphi$ is an antiautomorphism of $R$. Since $\chi$ commutes with $\wb{K}$, $\vphi$ preserves the $\wb{G}$-grading on $R$. Moreover, $\chi^2\in\Lambda$ implies that $\vphi^2$ acts trivially on the identity 
component of the $\wb{G}$-grading. Thus we can apply (for $\wb{G}$) the following result \cite[Proposition 6.4]{BZA}, whose proof does not require any assumptions about the characteristic:

\begin{proposition*}[Bahturin--Zaicev]
Let $R=M_n(F)$ be graded by a finite abelian group $G$. Let $\vphi$ be an antiautomorphism of $R$ that preserves the grading and acts as an involution on the identity component.
Then there exists an automorphism $\psi$ of $R$ that also preserves the grading such that $\vphi$ commutes with $\psi$ and $\vphi^2=\psi^2$.
\end{proposition*}

Now we can define a new $K$-action on $R$ by making $\chi$ act as $\psi$ (instead of $-\vphi$) and $\wb{K}$ as before. By construction, $R$ is a $K$-module algebra with 
respect to this new action, so the corresponding grading $R=\bigoplus_{g\in G}\wt{R}_g$ is a grading of $R$ as an associative algebra. Moreover, $*=\psi^{-1}\vphi$ is an involution on $R$ that preserves both gradings $R=\bigoplus_{g\in G}R_g$ and $R=\bigoplus_{g\in G}\wt{R}_g$. 

We need one auxiliary result that is a dualization of the so-called \emph{Exchange Theorem} of \cite{BG} and \cite{BSh} (which stems from \cite[Theorem 5.5]{BZA}). Incidentally, this dual form is valid without any restrictions on the base field and its proof is much simpler.

Suppose $R$ is a vector space, $G$ a group, and $R=\bigoplus_{g\in G}R_g$ and $R=\bigoplus_{g\in G}\wt{R}_g$ are two $G$-gradings. We will call these gradings {\em compatible} if 
for all $g\in G$, $\wt{R}_g=\bigoplus_{x\in G} (R_x\cap\wt{R}_g)$, or, equivalently, $R_g=\bigoplus_{x\in G} (\wt{R}_x\cap R_g)$. 

\begin{lemma}\label{exchange}
Let $R$ be a vector space with two compatible gradings $R=\bigoplus_{g\in G}R_g$ and $R=\bigoplus_{g\in G}\wt{R}_g$. Suppose $H\triangleleft G$ is such that the two factor-gradings 
by $G/H$ coincide. Set $R^h=\bigoplus_{g\in G}(\wt{R}_g\cap R_{gh})$. Then 
\[
R_g=\bigoplus_{h\in H}(\wt{R}_{gh^{-1}}\cap R^h).
\]
Moreover, if $R$ is a (nonassociative) algebra equipped with two such gradings and $H\subset Z(G)$, then $R=\bigoplus_{h\in H}R^h$ is an algebra grading.
\end{lemma}

\begin{proof}
Clearly,
\[
\wt{R}_{gh^{-1}}\cap R^h=\bigoplus_{x\in G} (\wt{R}_{gh^{-1}}\cap\wt{R}_x\cap R_{xh})=\wt{R}_{gh^{-1}}\cap R_{(gh^{-1})h}=\wt{R}_{gh^{-1}}\cap R_{g}.
\]
Thus
\begin{eqnarray*}
\bigoplus_{h\in H}(\wt{R}_{gh^{-1}}\cap R^h)&=&\bigoplus_{h\in H}(\wt{R}_{gh^{-1}}\cap R_{g})=\left(\bigoplus_{h\in H}\wt{R}_{gh^{-1}}\right)\cap R_{g}\\
&=&\wt{R}_{gH}\cap R_{g}=R_{gH}\cap R_{g}=R_g.
\end{eqnarray*}
Now if $R$ is an algebra graded in two ways and $H\subset Z(G)$, then for all $h_1, h_2\in H$ and $g_1,g_2\in G$, we have
\[
(\wt{R}_{g_1}\cap R_{g_1h_1})(\wt{R}_{g_2}\cap R_{g_2h_2})\subset \wt{R}_{g_1g_2}\cap R_{g_1h_1g_2h_2}=\wt{R}_{g_1g_2}\cap R_{(g_1g_2)h_1h_2}\subset R^{h_1h_2},
\]
which implies $R^{h_1}R^{h_2}\subset R^{h_1h_2}$.
\end{proof}

We apply Lemma \ref{exchange} in order to express $R_g$ in terms of 
$\wt{R}_g$ as follows. In our case $R^1=\bigoplus_{g\in G}(\wt{R}_g\cap R_g)=\bigoplus_{g\in G}\sks(\wt{R}_g,*)=\sks(R,*)$, 
$R^h=\bigoplus_{g\in G}(\wt{R}_g\cap R_{gh})=\bigoplus_{g\in G}\sym(\wt{R}_g,*)=\sym(R,*)$. Therefore,
\[
R_g=(\wt{R}_g\cap R^1)\oplus(\wt{R}_{gh}\cap R^h)=\sks(\wt{R}_g,*)\oplus\sym(\wt{R}_{gh},*).
\]
Restricting $R_g$ to $\rsl_n(F)$, we see that we have a grading of type II.

\end{proof}


\begin{thebibliography}{999}

\bibitem{BG}
Bahturin, Yuri; A. Giambruno. {\it Group Gradings on associative algebras with involution}, Canad. Math. Bull., to appear.
\bibitem{BSZ}
Bahturin, Yuri; S. Sehgal; M. Zaicev. {\it Group Gradings on Associative Algebras}, J. Algebra, {\bf 241} (2001), 677-698.
\bibitem{BSh}
Bahturin, Yuri; I. Shestakov. {\it Group gradings on associative superalgebras}, Contemporary Mathematics, \textbf{283}(2005), 849-868.
%\bibitem{BS}
%Y. Bahturin and I. Shestakov, {\it Gradings of Simple Jordan Algebras and Their Relation to the Gradings of Simple Associative algebras}, 
%Comm. Algebra, {\bf 999} (2002), 111-222.
\bibitem{BShZ}
Bahturin, Yuri; I. Shestakov; M. Zaicev. {\it Gradings on Simple Jordan and Lie  Algebras}, J. Algebra, {\bf 283} (2005), 849-868.
\bibitem{antaut}
Bahturin, Yuri; M. Zaicev. {\it Involutions on graded matrix algebras},  preprint.
\bibitem{surgrad}
Bahturin, Yuri;  M. Zaicev. {\it Graded algebras and graded identities},  Polynomial identities and combinatorial methods
(Pantelleria, 2001), 101-139, Lecture Notes in Pure and Appl. Math., 235, Dekker, New York, 2003.
\bibitem{BZnc}
Bahturin, Yuri;  M. Zaicev. {\it Group gradings on matrix algebras}, Canad. Math. Bull., {\bf 45} (2002), 499 - 508.
\bibitem{BZA}
Bahturin, Yuri; M. Zaicev. {\it Gradings on Simple Lie Algebras of Type ``A''}, J. Lie Theory, to appear.
\bibitem{BlMa}
Blau, P. S.; W. S. Martindale 3rd. {\it Lie isomorphisms in $*$-prime GPI rings with involution}, Taiwanese J. Math., {\bf 4} (2000), 215-252.
\bibitem{Dieu}
Dieudonn\'e, J. {\it Introduction to the theory of formal groups}. Pure and Applied Mathematics, 20. Marcel Dekker, Inc., New York, 1973.
\bibitem{Mart} 
Martindale 3rd., W. S. {\it Lie derivations of primitive rings}. Michigan Math. J., {\bf 11} (1964), 183--187. 
\bibitem{Mont}
Montgomery, S. {\it Hopf algebras and their actions on rings}. CBMS Regional Conference Series in Mathematics, 82. American Mathematical Society, Providence, RI, 1993.
%\bibitem{R} 
%L. H. Rowen, Polynomial Identities in Ring Theory, Academic Press, New York, 1980.
\end{thebibliography}
\end{document}